\documentclass[a4paper,12pt]{amsart}
\usepackage{ifthen}
\usepackage{mathrsfs}
\nonstopmode
\numberwithin{equation}{section}
\newtheorem{thm}{Theorem}[section]

\newtheorem{lem}[thm]{Lemma}
\newtheorem{prop}[thm]{Proposition}

\theoremstyle{definition}

\newenvironment{rem}{%
\bigskip
\noindent
\textsl{{\sl Remark. }}}{\bigskip}

\newenvironment{pf}[1][]{%
 \vskip 3mm
 \noindent
 \ifthenelse{\equal{#1}{}}%
  {{\slshape Proof. }}%
  {{\slshape #1.} }%
 }%
{\qed\bigskip}

\newcounter{alphabet}
\newcounter{tmp}
\newenvironment{Thm}[1][]{\refstepcounter{alphabet}%
\bigskip%
\noindent%
{\bf Theorem \Alph{alphabet}}%
\ifthenelse{\equal{#1}{}}{}{ (#1)}%
{\bf .}
\itshape}{\vskip 8pt}

\newcommand{\C}{{\mathbb C}}
\newcommand{\D}{{\mathbb D}}

\newcommand{\R}{{\mathbb R}}

\newcommand{\sphere}{{\widehat{\mathbb C}}}

\renewcommand{\Re}{{\operatorname{Re}\,}}

\newcommand{\Int}{{\operatorname{Int}\,}}

\newcommand{\inv}{^{-1}}

\renewcommand{\arg}{\,{\operatorname{arg}\,}}

\newcommand{\argl}{\,{\operatorname{arg}_\lambda\,}}
\newcommand{\Ula}{{U_{\lambda,\alpha}}}

\newcounter{minutes}
\divide\time by 60
\newcounter{hours}
\multiply\time by 60
\addtocounter{minutes}{-\time}

\begin{document}
\bibliographystyle{amsplain}
\title{
Quasiconformal extension of strongly spirallike functions
}


\author[T. Sugawa]{Toshiyuki Sugawa}
\address{Graduate School of Information Sciences,
Tohoku University, Aoba-ku, Sendai 980-8579, Japan}
\email{sugawa@math.is.tohoku.ac.jp}
\begin{abstract}
We show that a strongly $\lambda$-spirallike function of order $\alpha$
can be extended to a $\sin(\pi\alpha/2)$-quasiconformal automorphism
of the complex plane for $-\pi/2<\lambda<\pi/2$ and $0<\alpha<1$
with $|\lambda|<\pi\alpha/2.$
In order to prove it, we provide several geometric characterizations
of a strongly $\lambda$-spirallike domain of order $\alpha.$
We also give a concrete form of the mapping function of the standard
strongly $\lambda$-spirallike domain $U_{\lambda,\alpha}$ of order $\alpha.$
A key tool of the present study is the notion of $\lambda$-argument, which was
developed by Y.~C.~Kim and the author \cite{KS09spl}.
\end{abstract}
\thanks{
The author was supported in part by 
the JSPS Grant-in-Aid for Scientific Research (B), 22340025.
}
\maketitle

\section{Introduction and main result}

An analytic function $f$ on the unit disk $\D=\{z\in\C: |z|<1\}$ 
is called {\it starlike} if $f(0)=0$ and if $f$ maps $\D$ univalently
onto a domain starlike with respect to the origin.
It is well known that a non-constant analytic function $f$ on $\D$ with $f(0)=0$
is starlike precisely if $\Re zf'(z)/f(z)>0$
in $0<|z|<1.$
In what follows, the function $f(z)/z$ will be understood as an analytic
function $a_1+a_2z+\dots$ on $\D$ for an analytic function $f$ on $\D$ 
of the form $a_1z+a_2z^2+\dots$ on $\D$
so that the last condition makes sense for $z=0$ as well.

An analytic function $f$ on $\D$ is called 
{\it strongly starlike of order} $\alpha\in(0,1)$ if $f(0)=0,$
if $f$ is non-constant, and if
$$
\left|\arg\frac{zf'(z)}{f(z)}\right|<\frac{\pi\alpha}{2},
\quad z\in\D.
$$
The image $f(\D)$ of a strongly starlike function $f$ of order $\alpha$
is called a strongly starlike domain of order $\alpha$ (with respect to
the origin).
Various geometric characterizations of those domains are known.
See \cite{SugawaDual} for a summary of such characterizations.

Fait, Krzy\.z and Zygmunt \cite{FKZ76} gave quasiconformal extensions
of strongly starlike functions.

\begin{Thm}\label{Thm:FKZ}
A strongly starlike function of order $\alpha$ extends to a
$\sin(\pi\alpha/2)$-quasiconformal automorphism of $\C.$
\end{Thm}

Here, for $0\le k<1,$
a homeomorphism $f:\C\to\C$ is called $k$-quasiconformal if
$f$ has locally integrable partial derivatives on $\C$ (in the sense
of distributions) which satisfy $|\partial_{\bar z}f|\le k|\partial_zf|$ a.e.
Univalent functions on the unit disk with quasiconformal extension
are important in connection with the Bers embedding of Teichm\"uller
spaces. See \cite{Lehto:univ} for details.

Spirallike functions are natural generalization of starlike functions.
For a number $\lambda$ with $-\pi/2<\lambda<\pi/2,$ an analytic function
$f$ on $\D$ is called {\it $\lambda$-spirallike} (or {\it spirallike of
type $-\lambda$}, cf.~\cite{Pom:univ}) if $f(0)=0,$ if $f$ is univalent and
if $f(\D)$ is $\lambda$-spirallike (with respect to the origin).
Here, a domain $\Omega$ is said to be $\lambda$-spirallike (with respect
to the origin) if the $\lambda$-spiral segment
$$
[0,w]_\lambda:=\{w\exp(e^{i\lambda}t): t<0\}\cup\{0\}
$$
is contained in $\Omega$ whenever $w\in\Omega.$
Note that $0$-spirallikeness is nothing but starlikeness.
It is well known that a non-constant analytic function $f$ on $\D$
with $f(0)=0$ is $\lambda$-spirallike if and only if
$\Re e^{-i\lambda}zf'(z)/f(z)>0.$

Analogously, we can define the notion of strongly spirallike functions.
Let $\lambda\in(-\pi/2,\pi/2)$ and $\alpha\in(0,1).$
An analytic function $f$ on $\D$ will be called 
{\it strongly $\lambda$-spirallike of order $\alpha$} if $f(0)=0,$
if $f$ is non-constant, and if
\begin{equation}\label{eq:sspl}
\left|\arg\frac{zf'(z)}{f(z)}-\lambda\right|<\frac{\pi\alpha}{2},
\quad z\in\D.
\end{equation}
Since $zf'(z)/f(z)$ assumes the value $1$ at $z=0,$ we have the
inevitable constraint $|\lambda|<\pi\alpha/2.$
The image $f(\D)$ of $\D$ under a strongly $\lambda$-spirallike
function $f$ of order $\alpha$ 
is called a strongly $\lambda$-spirallike domain of order $\alpha$
(with respect to the origin).
If we do not need to specify $\lambda,$ we will simply call it
strongly spirallike of order $\alpha.$
Our principal aim in the present note is to extend Theorem A
for strongly spirallike functions.

\begin{thm}\label{thm:main}
A strongly spirallike function $f$ of order $\alpha\in(0,1)$
extends to a $\sin(\pi\alpha/2)$-quasiconformal 
automorphism of $\C.$
\end{thm}

\begin{rem}
Betker \cite{Bet92} proved the following: {\sl If a L\"owner chain $f(z,t)$
satisfies $\partial_tf(z,t)=z\partial_zf(z,t)p(z,t)$ 
with $|\arg p(z,t)|\le \pi\alpha/2, z\in\D, t\ge0,$
then $f(z)=f(z, 0)$ extends to a $\sin(\pi\alpha/2)$-quasiconformal
automorphism of $\sphere.$}

Therefore, if we could find such a L\"owner chain for a strongly 
$\lambda$-spirallike function $f$ of order $\alpha,$
Theorem \ref{thm:main} would follow from Betker's theorem
(by replacing $\C$ by $\sphere$).
The author, however, is not able to find
such a L\"owner chain for general $\lambda$ so far.
For instance, a standard L\"owner chain for a $\lambda$-spirallike function $f$
is given by
$$
f(z,t)=e^{(1+i\tan\lambda)t}f(e^{-it\tan\lambda}z), \quad z\in\D, t\ge0
$$
(see \cite[Theorem 6.6]{Pom:univ}).
Since
$$
p(z,t)=\frac{\partial_tf(z,t)}{z\partial_zf(z,t)}
=\frac{e^t}{\cos\lambda}\left[
\frac{f(z)}{e^{-i\lambda}zf'(z)}+i\sin\lambda
\right],
$$
this chain does not necessarily satisfy the assumption of Betker's theorem for
a strongly $\lambda$-spirallike function $f$ of order $\alpha$
unless $\lambda=0.$
\end{rem}

The proof of our theorem will be given along the line same as in \cite{FKZ76}.
We shall use, however, the idea of $\lambda$-argument developed in \cite{KS09spl}
for a substitute of the usual angle.

In Section 2, we give a couple of geometric characterizations
of strongly $\lambda$-spirallike domains of order $\alpha.$
By using it, we prove Theorem \ref{thm:main} in Section 3.
We also give additional remarks in Section 3.

\section{Characterizations of strongly spirallike functions}

Let $\lambda$ be a real number with $|\lambda|<\pi/2.$
The $\lambda$-argument of a complex number $w\ne0$ is defined to be
$\arg w-\tan\lambda\cdot\log|w|$ and is denoted by $\argl w.$
In other words, $\theta=\argl w$ if and only if
$w$ lies on the $\lambda$-logarithmic spiral ($\lambda$-spiral for short)
$e^{i\theta}\gamma_\lambda,$
here $\gamma_\lambda=\{\exp{e^{i\lambda}t}: t\in\R\}.$
We denote by $P_\lambda(r,\theta)$ the complex number $w$ with
$r=|w|$ and $\theta=\argl w;$ namely,
$$
P_\lambda(r,\theta)=re^{i(\theta+\tan\lambda\log r)}
=\exp\big[(1+i\tan\lambda)\log r+i\theta\big].
$$
Here we define $P_\lambda(0,\theta)=0$ for any $\theta.$
The following elementary property is sometimes useful:
\begin{equation}\label{eq:P}
P_\lambda(r,\theta)P_\lambda(s,t)=P_\lambda(rs,\theta+t).
\end{equation}

For a domain $\Omega\subset\C$ with $0\in\Omega,$ we define a periodic
function $R_\lambda:\R\to(0,+\infty]$ with period $2\pi$ by
$$
R_\lambda(\theta)=\sup\{r>0: [0,P_\lambda(r,\theta)]_\lambda\subset\Omega\},
$$
which will be called the radius function of $\Omega$ with respect to
the $\lambda$-spirals.
We notice that $R_\lambda$ is lower semi-continuous, namely,
$\liminf_{t\to\theta}R_\lambda(t)\ge R_\lambda(\theta).$
Conversely, if $R:\R\to(0,+\infty]$ is a lower semi-continuous periodic
function with period $2\pi,$ then
$\Omega=\{P_\lambda(r,\theta): \theta\in\R, 0\le r<R(\theta)\}$
is a $\lambda$-spirallike domain with $0\in\Omega.$

If $R_\lambda$ is continuous at every point, then
$\Omega$ must be $\lambda$-spirallike.
This was stated in \cite{SugawaDual} without proof for the case
when $\lambda=0.$
We supply a proof for it in the present note.

\begin{lem}\label{lem:cont}
Let $\Omega$ be a plane domain with $0\in\Omega.$
If $R_\lambda:\R\to(0,+\infty]$ is continuous, then
$\Omega$ is $\lambda$-spirallike with respect to the origin.
If furthermore $R_\lambda(\theta)$ is finite for every $\theta,$
$\Omega$ is a bounded Jordan domain.
\end{lem}

\begin{pf}
Let $\Omega_0$ be the subdomain of $\Omega$ given as
$\{P_\lambda(r,\theta): \theta\in\R, 0\le r<R_\lambda(\theta)\}.$
It is enough to show that $\Omega_0=\Omega.$
Suppose, to the contrary, that $\Omega_0\ne\Omega.$
Then there is a point $w_0\in\partial\Omega_0\cap\Omega.$
By the definition of $R_\lambda,$ we must have $|w_0|>R_\lambda(\theta_0)$
for $\theta_0=\argl w_0.$
Since $w_0\in\partial\Omega_0,$there is a sequence $w_n$ in $\Omega_0$
converging to $w_0.$
Let $\theta_n=\argl w_n$ so that $\theta_n\to\theta_0.$
Then 
$$
\limsup_{n\to\infty}R_\lambda(\theta_n)
\ge \lim_{n\to\infty}|w_n|= |w_0|>R_\lambda(\theta_0),
$$
which violates continuity of $R_\lambda.$
Thus we have shown that $\Omega$ is $\lambda$-spirallike.
If $R_\lambda$ is finite, the correspondence
$e^{i\theta}\mapsto P_\lambda(R_\lambda(\theta),\theta)$ gives
an injective continuous mapping from $S^1=\partial\D$ into $\C.$
Therefore, we conclude that $\Omega$ is a bounded Jordan domain.
\end{pf}

We defined strong $\lambda$-spirallikeness of order $\alpha$ for domains
in terms of mapping functions.
It is possible to give geometric characterizations of such domains.
To state it, we need the standard $\lambda$-spirallike domain
\begin{align*}
\Ula&=\{P_\lambda(r,\theta): 0\le\theta<2\pi,~0\le r
<\max\{e^{B\theta}, e^{-A(2\pi-\theta)}\}\} \\
&=\{\exp[(1+i\tan\lambda)t+i\theta]: 0\le\theta<2\pi,~
t<\max\{B\theta,-A(2\pi-\theta)\}\}\cup\{0\}, 
\end{align*}
where $A$ and $B$ are defined by
\begin{equation}\label{eq:AB}
\begin{cases}
A&=-\sin\lambda\cos\lambda+\cos^2\lambda\tan\tfrac{\pi\alpha}2
=\cos^2\lambda\left(\tan\tfrac{\pi\alpha}2-\tan\lambda\right), \\
B&=-\sin\lambda\cos\lambda-\cos^2\lambda\tan\tfrac{\pi\alpha}2
=-\cos^2\lambda\left(\tan\tfrac{\pi\alpha}2+\tan\lambda\right).
\end{cases}
\end{equation}
Note that $A>0$ and $B<0.$
We observe that $B\theta\ge -A(2\pi-\theta)$ if and only if $\theta\le\theta^*$
for $0\le\theta\le2\pi,$ where
\begin{equation}\label{eq:theta}
\theta^*=\frac{2\pi A}{A-B}=\pi\left(1-\frac{\tan\lambda}{\tan(\pi\alpha/2)}
\right).
\end{equation}
We shall use the following parametrization of the boundary of $\Ula:$
\begin{equation}\label{eq:para}
\beta(\theta)=\begin{cases}
e^{[i+B(1+i\tan\lambda)]\theta} & \quad (0\le\theta\le\theta^*) \\
e^{[i+A(1+i\tan\lambda)]\theta} &\quad (\theta^*-2\pi<\theta<0).
\end{cases}
\end{equation}
In view of \eqref{eq:P}, for $w_0=P_\lambda(r_0,\theta_0),$ we observe
$$
w_0\Ula=\{P_\lambda(r,\theta): \theta_0\le\theta<\theta_0+2\pi,~
0\le r<r_0\max\{e^{B(\theta-\theta_0)}, e^{-A(2\pi-\theta+\theta_0)}\}\}.
$$
In particular, $w_0\Ula\subset w_1\Ula$ 
whenever $|w_0|<|w_1|$ and $\argl w_0=\argl w_1.$

\begin{thm}\label{thm:char}
Let $\lambda$ and $\alpha$ be real numbers with $|\lambda|<\pi\alpha/2<\pi/2.$
For a domain $\Omega$ in $\C$ with $0\in\Omega,$
each condition in the following implies the others.
\begin{enumerate}
\item[(i)]
There exists a strongly $\lambda$-spirallike function
$f$ of order $\alpha$ such that $\Omega=f(\D).$
\item[(ii)]
The radius function $R_\lambda$ of $\Omega$ with respect to $\lambda$-spirals
is bounded, absolutely continuous on $[0,2\pi]$ and satisfies
\begin{equation}\label{eq:double}
\left|\frac{R_\lambda'(\theta)}{R_\lambda(\theta)}+\sin\lambda\cos\lambda
\right|\le \cos^2\lambda\tan\frac{\pi\alpha}2
\end{equation}
for almost every $\theta.$
\item[(iii)]
$wU_{\lambda,\alpha}\subset\Omega$ whenever $w\in\partial\Omega.$
\item[(iv)]
$wU_{\lambda,\alpha}\subset\Omega$ whenever $w\in\Omega.$
\end{enumerate}
\end{thm}

For the proof, we first need the following result \cite{KS11PT}.
Of course, this can also be shown directly (see \cite{BK69}).

\begin{lem}\label{lem:bdd}
A strongly spirallike function of order $\alpha<1$ is bounded.
\end{lem}

We also need results in the theory of Hardy spaces.
We extract necessary information from \S\S 3.4-5 of \cite{Duren:hp}.

\begin{lem}\label{lem:Privalov}
Suppose that $f'\in H^1$ for an analytic function $f$ on the unit disk $\D.$
Then the boundary function $f(e^{it})$ of $f$ is absolutely continuous
and satisfies $\frac{d}{dt}f(e^{it})=ie^{it}f'(e^{it})\ne0$ a.e., where
$f'(e^{it})$ is the nontangential limit of $f'.$
In particular, for a Borel subset $E$ of $\partial\D,$
$f(E)$ is of linear measure zero if and only if so is $E.$
\end{lem}

\begin{pf}[Proof of Theorem \ref{thm:char}]
Basically, we can modify arguments in \cite{SugawaDual} to prove
the theorem.

We first show that (i) implies (ii).
Assume that $\Omega=f(\D)$ for a strongly $\lambda$-spirallike
function $f$ of order $\alpha.$
Then $q(z)=zf'(z)/f(z)$ is subordinate to the function
\begin{equation}\label{eq:Q}
Q(z)=\left(\frac{1+ze^{2i\lambda/\alpha}}{1-z}\right)^\alpha.
\end{equation}
Since $Q$ belongs to the Hardy space $H^p$ for $p<1/\alpha,$
we have $q\in H^p$ for $p<1/\alpha$ by Littlewood's subordination theorem
\cite[Theorem 1.7]{Duren:hp}.
In particular, $q\in H^1.$
Since $f(z)/z$ is bounded by Lemma \ref{lem:bdd},
$f'=qf/z$ belongs to $H^1,$ too.
Thus by the above lemma,
$f$ has an absolutely continuous boundary function $f(e^{it})$
and $f'$ has nontangential limit $f'(e^{it})$ which vanishes
almost nowhere.
In particular, the boundary function of $f$ can be described by
$f(e^{it})=\rho(t)e^{i\Theta(t)}$ with absolutely continuous
functions $\rho:\R\to(0,+\infty)$ and $\Theta:\R\to\R$
with $\Theta(t+2\pi)\equiv\Theta(t)+2\pi.$
Let $\Theta_\lambda(t)$ be the $\lambda$-argument of $f(e^{it});$
more explicitly, it can be given by
$$
\Theta_\lambda(t)=\Theta(t)-\tan\lambda\cdot\log\rho(t).
$$
Then
$$
q(e^{it})=\frac{e^{it}f'(e^{it})}{f(e^{it})}
=\Theta'(t)-i\frac{\rho'(t)}{\rho(t)}
=\Theta_\lambda'(t)-i(1+i\tan\lambda)\frac{\rho'(t)}{\rho(t)}
$$
and thus
$$
e^{-i\lambda}q(e^{it})
=e^{-i\lambda}\Theta_\lambda'(t)-\frac{i\rho'(t)}{\rho(t)\cos\lambda}.
$$
The last quantity lies in the closed sector $|\arg w|\le \pi\alpha/2$
for almost every $t$ by \eqref{eq:sspl}.
Since $q(e^{it})=e^{it}f'(e^{it})/f(e^{it})\ne0$ a.e., 
one can conclude that $\Theta_\lambda'(t)>0$ for a.e.~$t.$
In particular, $\Theta_\lambda:\R\to\R$ is a homeomorphism
and $\Theta_\lambda\inv$ preserves the sets of linear measure zero;
in other words, $\Theta_\lambda\inv$ is absolutely continuous.
By definition, $R_\lambda$ can be expressed by 
$R_\lambda(\theta)=\rho(\Theta_\lambda\inv(\theta)).$
Therefore, $R_\lambda$ is also absolutely continuous and
satisfies the relation 
$R_\lambda'(\Theta_\lambda(t))=\rho'(t)/\Theta_\lambda'(t)$
for a.e.~$t.$
Hence,
$$
e^{-i\lambda}q(e^{it})=\frac{\Theta_\lambda'(t)}{\cos\lambda}
\left[\cos^2\lambda-i\left(\frac{R_\lambda'(\theta)}{R_\lambda(\theta)}
+\sin\lambda\cos\lambda\right)\right]
$$
for $\theta=\Theta_\lambda(t).$
Since the argument of $e^{-i\lambda}q$ is between $-\pi\alpha/2$
and $\pi\alpha/2,$ we now obtain \eqref{eq:double} for a.e.~$\theta.$

Secondly, we show that (ii) implies (iii).
Assume (ii).
Then $\Omega$ is $\lambda$-spirallike by Lemma \ref{lem:cont}.
By integrating the inequality in \eqref{eq:double}, we have
$$
B(\theta_2-\theta_1)\le \log R_\lambda(\theta_2)-\log R_\lambda(\theta_1)
\le A(\theta_2-\theta_1)
$$
for $\theta_1<\theta_2,$ where $A$ and $B$ are the constants given
in \eqref{eq:AB}.
Let $\theta_0\in\R$ and $0<\theta<2\pi.$
Applying the above inequalities to combinations of the three numbers
$\theta_0+\theta-2\pi<\theta_0<\theta_0+\theta,$ we obtain
$$
\max\{e^{B\theta}, e^{-A(2\pi-\theta)}\}
\le\frac{R_\lambda(\theta_0+\theta)}{R_\lambda(\theta_0)},
$$
which implies that $w_0\Ula\subset\Omega$
for $w_0=P_\lambda(R_\lambda(\theta_0),\theta_0).$
Thus condition (iii) has been checked.

Thirdly, we show that (iii) implies (iv).
For a point $w=P_\lambda(r,\theta_0)\in\Omega,$
let $w_0=P_\lambda(R_\lambda(\theta_0),\theta_0).$
Then $w\Ula\subset w_0\Ula\subset\Omega$
(recall the observation made right before the statement of 
Theorem \ref{thm:char}).
We thus have shown (iii) $\Rightarrow$ (iv).

Finally, we show that (iv) implies (i) by
following the argument of Ma and Minda \cite{MaMinda91}.
Let us assume (iv).
Since $w\Ula\cup\{w\}$ contains the $\lambda$-spiral segment $[0,w]_\lambda$
for $w\in\Omega,$
it implies that $\Omega$ is $\lambda$-spirallike and, in particular,
simply connected.
Now we can take a conformal homeomorphism $f:\D\to\Omega$
with $f(0)=0, f'(0)>0.$
We may assume that $f'(0)=1.$

First we show that $f(z)\Ula\subset f(\D_r),$
where $r=|z|$ and $\D_r=\{z: |z|<r\}.$
Indeed, for a fixed $w\in\Ula,$ we have
$f(z)w\in f(z)\Ula\subset\Omega$ by (iv).
Thus $g(z)=f\inv(f(z)w)$ is an analytic function in $|z|<1$
and satisfies $|g(z)|<1$ and $g'(0)=w\in\Ula\subset\D.$
Therefore, by Schwarz's lemma, we obtain $|g(z)|<|z|$
for $0<r=|z|<1,$ which implies $f(z)w=f(g(z))\in f(\D_r)$
as required.

Fix an arbitrary point $z_0\in\D$ with $z_0\ne0$
and let $\gamma(t)=f(z_0e^{it})$ and $w_0=f(z_0).$
By the above observation, we see that the curve $\gamma$
encloses the domain $w_0\Ula.$ 
Thus, we have the inequalities
$\arg w_0\beta'(0-)\le\arg\gamma'(0)\le\arg w_0\beta'(0+);$
in other words, $\arg [w_0\beta'(0+)/\gamma'(0)]\ge0$
and $\arg [w_0\beta'(0-)/\gamma'(0)]\le 0,$
where $\beta$ is given in \eqref{eq:para}.
Letting $q(z)=zf'(z)/f(z)$ as before,
we now compute
\begin{align*}
\frac{w_0\beta'(0+)}{\gamma'(0)}
&=\frac{i+B(1+i\tan\lambda)}{iq(z_0)}
=\frac{e^{-i\lambda}-iB/\cos\lambda}{e^{-i\lambda}q(z_0)} \\
&=\frac{\cos\lambda(1+i\tan(\pi\alpha/2))}{e^{-i\lambda}q(z_0)}
=\frac{\cos\lambda}{\cos(\pi\alpha/2)}\cdot
\frac{e^{i\pi\alpha/2}}{e^{-i\lambda}q(z_0)}
\end{align*}
Therefore,
$$
\arg\frac{w_0\beta'(0+)}{\gamma'(0)}=-\arg q(z_0)+\lambda+\frac{\pi\alpha}2
\ge0.
$$
Similarly,
$$
\arg\frac{w_0\beta'(0-)}{\gamma'(0)}
=\frac{i+A(1+i\tan\lambda)}{iq(z_0)}
=-\arg q(z_0)+\lambda-\frac{\pi\alpha}2\le0.
$$
Therefore, we have $-\pi\alpha/2\le\arg q(z_0)-\lambda\le\pi\alpha/2.$
Since $\arg q$ is harmonic on $\D,$ the inequalities are both strict.
In other words, $f$ is strongly $\lambda$-spirallike
of order $\alpha.$
In this way, we have verified condition (i).
\end{pf}

\section{Quasiconformal extension and additional observations}

In order to prove Theorem \ref{thm:main}, it is enough to construct
a $\sin(\pi\alpha/2)$-quasiconformal reflection $h$ with $h(0)=\infty$
in the boundary of a strongly spirallike domain of order $\alpha.$
Here, a self-mapping $h$ of the Riemann sphere 
$\sphere=\C\cup\{\infty\}$ is called a $k$-quasiconformal reflection
in the Jordan curve $C$ if $h$ is an orientation-reversing homeomorphism
such that $h(z)=z$ for every $z\in C,$ that $h=h\inv$
and that $h$ has locally integrable
partial (distributional) derivatives on $\C\setminus\{h(\infty)\}$
with $|\partial_zh|\le k|\partial_{\bar z}h|$ a.e.
Indeed, a conformal mapping $f$ of the unit disk $\D$
onto the interior $\Omega$ of $C$ can be extended to $k$-quasiconformal
mapping $\tilde f$ of $\sphere$ by setting
$$
\tilde f(z)=\begin{cases}
f(z) &\quad (z\in\overline\Omega) \\
h(f(1/\bar z)) &\quad (z\in\sphere\setminus\Omega)
\end{cases}
$$
(see \cite{Lehto:univ} for details).
Note that the above $\tilde f$ satisfies $\tilde f(\infty)=\infty$
when $f(0)=0$ and $h(0)=\infty.$

Let $\Omega$ be a strongly $\lambda$-spirallike domain of order
$\alpha$ with $|\lambda|<\pi\alpha/2<\pi/2$
and let $R_\lambda$ be the radius function of $\Omega$ with respect
to $\lambda$-spirals.
Then, by Theorem \ref{thm:char}, $R_\lambda$ is a bounded, absolutely
continuous function satisfying \eqref{eq:double}.
Therefore, Lemma \ref{lem:cont} implies that $\Omega$ is a bounded Jordan
domain.
We now define a reflection $h$ in $\partial\Omega$ by
$$
h(P_\lambda(r,\theta))=P_\lambda(R_\lambda(\theta)^2/r, \theta).
$$
Here, we interpret that $h(0)=\infty, h(\infty)=0$ when $r=0,+\infty,$
respectively.
It is easy to see that $h$ has locally integrable partial derivatives
(in the sense of distributions) on $\C\setminus\{0\}.$
To estimate the dilatation of $w=h(z),$ we use the logarithmic coordinates
$Z=X+iY=\log z$ and $W=U+iV=\log w.$
By the conformal invariance of the dilatation, we note that
$|\partial_Z W/\partial_{\overline Z}W|
=|\partial_z w/\partial_{\bar z}w|.$
For short, we set $\rho(\theta)=\log R_\lambda(\theta).$
Using the relation $\theta=\argl z=Y-X\tan\lambda,$ 
we can express $W=\log h$ by
$$
W=2\rho(Y-X\tan\lambda)-X
+i\left[(Y-X\tan\lambda)+(2\rho(Y-X\tan\lambda)-X)\tan\lambda\right].
$$
Therefore,
\begin{align*}
\partial_X W&
=-2(1+i\tan\lambda)\rho'\tan\lambda-(1+2i\tan\lambda), \\
\partial_Y W&
=2(1+i\tan\lambda)\rho'+i,
\end{align*}
and, after simplifications, one has
$$
\frac{\partial_Z W}{\partial_{\overline Z}W}
=\frac{\partial_XW-i\,\partial_YW}{\partial_XW+i\,\partial_YW}
=-e^{-2i\lambda}\frac{\rho'+\sin\lambda\cos\lambda}%
{\rho'+\sin\lambda\cos\lambda+i\cos^2\lambda}.
$$
Since $|\rho'+\sin\lambda\cos\lambda|\le \cos^2\lambda\tan(\pi\alpha/2)$
by \eqref{eq:double}, we finally obtain
$$
\left|\frac{\partial_Z W}{\partial_{\overline Z}W}\right|
\le \frac{\cos^2\lambda\tan(\pi\alpha/2)}%
{\sqrt{[\cos^2\lambda\tan(\pi\alpha/2)]^2+\cos^4\lambda}}
=\sin\frac{\pi\alpha}2,
$$
which implies that $h$ is a $\sin(\pi\alpha/2)$-quasiconformal
reflection in $\partial\Omega.$
We have finished the proof of Theorem \ref{thm:main}.

\bigskip

We end the present note by giving additional observations.
In \cite{SugawaDual}, the author emphasized a self-duality of strong
starlikeness.
For a bounded domain $\Omega$ in $\C$ with $0\in\Omega$ we set
$\Omega^\vee=\{w: 1/w\in\sphere\setminus\overline\Omega\}.$
Note that $(\Omega^\vee)^\vee=\Omega$ when $\Int\overline\Omega=\Omega.$
The self-duality means that $\Omega$ is strongly starlike of order $\alpha$
if and only if so is $\Omega^\vee.$
One can extend this idea to strong spirallikeness.
Suppose that the radius function $R_\lambda$ of a domain $\Omega$
with $0\in\Omega$ is continuous and bounded.
Then $\Omega$ is a $\lambda$-spirallike Jordan domain by Lemma \ref{lem:cont} and
$\Omega=\{P_\lambda(r,\theta): \theta\in\R,~r<R_\lambda(\theta)\}.$
By the relation $1/P_\lambda(r,\theta)=P_\lambda(1/r, -\theta)$
which follows from \eqref{eq:P}, the domain $\Omega^\vee$
can be described by 
$\Omega^\vee=\{P_\lambda(r,\theta): r<1/R_\lambda(-\theta)\}.$
Therefore, $\Omega^\vee$ is $\lambda$-spirallike and
the radius function $R_\lambda^\vee$ of $\Omega^\vee$
is given by $R_\lambda^\vee(\theta)=1/R_\lambda(-\theta).$
In particular, $(\log R_\lambda^\vee)'(\theta)=(\log R_\lambda)'(-\theta)$
at every differentiable point.
Theorem \ref{thm:char} now yields the following.

\begin{thm}\label{thm:dual}
Let $\Omega\subset\C$ be a domain with $0\in\C.$
If $\Omega$ is a strongly $\lambda$-spirallike domain of order $\alpha,$
then so is $\Omega^\vee=\{w: 1/w\in\sphere\setminus\overline\Omega\}.$
\end{thm}

The radius function $R_\lambda$ of
the standard domain $U_{\lambda,\alpha}$ is absolutely continuous
and satisfies the equation
$$
(\log R_\lambda)'(\theta)=\begin{cases}
B=-\sin\lambda\cos\lambda-\cos^2\lambda\tan(\pi\alpha/2)
&\quad (0<\theta<\theta^*) \\
A=-\sin\lambda\cos\lambda+\cos^2\lambda\tan(\pi\alpha/2)
&\quad (\theta^*-2\pi<\theta<0).
\end{cases}
$$
Therefore, by Theorem \ref{thm:char}, $U_{\lambda,\alpha}$ itself is
a strongly $\lambda$-spirallike domain of order $\alpha.$

In the context of duality, it is natural to consider the strongly
$\lambda$-spirallike domain $V_{\lambda,\alpha}=U_{\lambda,\alpha}^\vee$
of order $\alpha.$
Explicitly,
$$
V_{\lambda,\alpha}=\{P_\lambda(r,\theta): 
0\le\theta<2\pi, r<\min\{e^{A\theta},e^{B(\theta-2\pi)}\}\}.
$$
Its radius function $R_\lambda^\vee(\theta)=1/R_\lambda(-\theta)$
satisfies
$$
(\log R_\lambda^\vee)'(\theta)=\begin{cases}
A&\quad (0<\theta<2\pi-\theta^*) \\
B&\quad (-\theta^*<\theta<0).
\end{cases}
$$
We note that $V_{\lambda,\alpha}$ is similar to $U_{\lambda,\alpha}.$
Indeed, we have $U_{\lambda,\alpha}=w^*V_{\lambda,\alpha},$
where $w^*=P_\lambda(e^{B\theta^*},\theta^*)$ is the other
tip of $U_{\lambda,\alpha}$ than $1.$
We now have more characterizations of strongly $\lambda$-spirallike
domains of order $\alpha.$

\begin{thm}
Let $\lambda$ and $\alpha$ be real numbers with $|\lambda|<\pi\alpha/2
<\pi/2.$
For a plane domain $\Omega$ with $0\in\Omega$ and $\Omega\ne\C,$
the following conditions are equivalent:
\begin{enumerate}
\item[(i)]
$\Omega$ is strongly $\lambda$-spirallike of order $\alpha.$
\item[(ii)]
$\Omega\subset wV_{\lambda,\alpha}$ whenever $w\in\partial\Omega.$
\item[(iii)]
$\Omega\subset wV_{\lambda,\alpha}$ whenever $w\in\C\setminus
\overline\Omega.$
\end{enumerate}
\end{thm}

\begin{pf}
First assume (i).
Then $\Omega^\vee$ is strongly $\lambda$-spirallike of order $\alpha$
by Theorem \ref{thm:dual}.
Theorem \ref{thm:char} now implies that $wU_{\lambda,\alpha}\subset
\Omega^\vee$ for $w\in\Omega^\vee.$
Thus $\Omega\subset (wU_{\lambda,\alpha})^\vee=w'V_{\lambda,\alpha}$
for $w'=1/w\in\C\setminus\overline\Omega.$
In this way, (iii) follows from (i).
Obviously, (iii) implies (ii).

Finally, assume (ii).
Then it is easy to see that $\Omega$ is a $\lambda$-spirallike bounded
domain with the property $\Int\overline\Omega=\Omega.$
Thus (ii) is equivalent to the condition that $(1/w)U_{\lambda,\alpha}
\subset\Omega^\vee$ for $w\in\partial\Omega.$
Since the correspondence $w\mapsto 1/w$ gives a bijection
from $\partial\Omega$ onto $\partial\Omega^\vee,$
Theorem \ref{thm:char} implies that $\Omega^\vee$ is strongly
$\lambda$-spirallike of order $\alpha.$
Hence, so is $\Omega=(\Omega^\vee)^\vee$ by Theorem \ref{thm:dual}.
We have shown that (ii) implies (i).
\end{pf}

As in the case of strongly starlike functions, the
function $k_{\lambda,\alpha}$ on $\D$ determined by the differential equation
$zk_{\lambda,\alpha}'(z)/k_{\lambda,\alpha}(z)=Q(z)$
and the initial conditions $k_{\lambda,\alpha}(0)=0, k_{\lambda,\alpha}'(0)=1$
is expected to be extremal in many respects among the class
of strongly $\lambda$-spirallike functions of order $\alpha.$
Here, $Q$ is given in \eqref{eq:Q} and therefore $k_{\lambda,\alpha}$
can be expressed explicitly by
$$
k_{\lambda,\alpha}(z)
=z\exp\int_0^z
\left[\left(\frac{1+\zeta e^{2i\lambda/\alpha}}{1-\zeta}\right)^\alpha-1\right]
\frac{d\zeta}{\zeta}.
$$
We will see that the image domain $k_{\lambda,\alpha}(\D)$
is a dilation of $U_{\lambda,\alpha}$ about the origin.
Indeed, $g=k_{\lambda,\alpha}/k_{\lambda,\alpha}(1)$ gives
a conformal representation of $U_{\lambda,\alpha}$ such that
$g(0)=0$ and $g(1)=1$ as we assert in the next theorem.

\begin{thm}\label{thm:U}
The function
$$
g(z)=z\exp\int_1^z
\left[\left(\frac{1+\zeta e^{2i\lambda/\alpha}}{1-\zeta}\right)^\alpha-1\right]
\frac{d\zeta}{\zeta}
$$
maps $\D$ conformally onto $U_{\lambda,\alpha}$
in such a way that $g(0)=0$ and $g(1)=1$
for $\lambda$ and $\alpha$ with $|\lambda|<\pi\alpha/2<\pi/2.$
\end{thm}

\begin{pf}
We write $k=k_{\lambda,\alpha}$ for brevity.
It is enough to see that $k(\D)=w_0U_{\lambda,\alpha}$
for the point $w_0=k(1).$
By definition, $k$ is strongly $\lambda$-spirallike of order $\alpha.$
We denote by $R_\lambda$ the radius function of $k(\D)$
with respect to $\lambda$-spirals.
Then, $R_\lambda$ is absolutely continuous and bounded by 
Theorem \ref{thm:char}.

We first recall that the function $Q$ given in \eqref{eq:Q} maps
$\D$ onto the sector $|\arg w-\lambda|<\pi\alpha/2.$
By the form of $Q,$ we observe that $\arg Q(e^{it})=\lambda+\pi\alpha/2$
for $0<t<t_1$ and that $\arg Q(e^{it})=\lambda-\pi\alpha/2$
for $t_1-2\pi<t<0,$ where $t_1=\pi-2\lambda/\alpha.$
If we write $k_{\lambda,\alpha}(e^{it})
=\rho(t)\exp i(\Theta_\lambda(t)+\tan\lambda\cdot\log\rho(t))$
as in the proof of Theorem \ref{thm:char}, we have the relation
$$
e^{-i\lambda}Q(e^{it})=\frac{\Theta_\lambda'(t)}{\cos\lambda}
\left[\cos^2\lambda-i\left(\frac{R_\lambda'(\theta)}{R_\lambda(\theta)}
+\sin\lambda\cos\lambda\right)\right]
$$
for $\theta=\Theta_\lambda(t).$
We set $\theta_0=\Theta_\lambda(0)$ and $\theta_1=\Theta_\lambda(t_1).$
Then, the above observation yields the relation
$$
\frac{R_\lambda'(\theta)}{R_\lambda(\theta)}
=\begin{cases}
B &\quad (\theta_0<\theta<\theta_1), \\
A &\quad (\theta_1-2\pi<\theta<\theta_0),
\end{cases}
$$
where $A$ and $B$ are given in \eqref{eq:AB}.
By integrating it, we have 
$$
R_\lambda(\theta)=
\begin{cases}
R_\lambda(\theta_0)e^{B(\theta-\theta_0)}
&\quad (\theta_0<\theta<\theta_1), \\
R_\lambda(\theta_0)e^{A(\theta-\theta_0)}
&\quad (\theta_1-2\pi<\theta<\theta_0).
\end{cases}
$$
Since $R_\lambda(\theta_1)=R_\lambda(\theta_1-2\pi),$
we have $e^{B(\theta_1-\theta_0)}=e^{A(\theta_1-\theta_0-2\pi)},$
equivalently, $\theta_1-\theta_0=\theta^*,$ where
$\theta^*$ is given in \eqref{eq:theta}.
We have thus seen that $R_\lambda(\theta+\theta_0)/R_\lambda(\theta_0)$
is the same radius function of the standard domain $U_{\lambda,\alpha}$
with respect to $\lambda$-spirals.
Hence, we conclude that $k(\D)=w_0U_{\lambda,\alpha}$
for the point $w_0=k(1).$
\end{pf}

As a by-product of the above proof, we can evaluate
a definite integral.
With the notation in the proof, we see that $g(e^{it_1})$
is the other tip $w^*=P_\lambda(e^{B\theta^*},\theta^*)$
of $U_{\lambda,\alpha}.$
Namely, $g(e^{it_1})=\exp[i+B(1+i\tan\lambda)]\theta^*.$
In view of the expression of $g,$ we obtain the relation
$$
it_1+\int_1^{e^{it_1}}
\left[\left(\frac{1+\zeta e^{2i\lambda/\alpha}}{1-\zeta}\right)^\alpha-1\right]
\frac{d\zeta}{\zeta}
=[i+B(1+i\tan\lambda)]\theta^*.
$$
We take $\zeta=e^{it}~(0\le t\le t_1)$ 
as the path of integration to obtain
$$
ie^{i(\lambda+\pi\alpha/2)}\int_0^{t_1}
\left(\frac{\cos(\frac t2+\frac\lambda\alpha)}{\sin\frac t2}\right)^\alpha
dt=[i+B(1+i\tan\lambda)]\theta^*.
$$
Substituting $t_1=\pi-2\lambda/\alpha,$ we obtain
$$
\int_0^{\pi-2\lambda/\alpha}
\frac{\cos^\alpha(\frac t2+\frac\lambda\alpha)}{\sin^\alpha\frac t2}dt
=e^{-i(\lambda+\pi\alpha/2)}[1-iB(1+i\tan\lambda)]\theta^*
=\frac{2\pi\sin(\frac{\pi\alpha}2-\lambda)}{\sin\pi\alpha}.
$$
By the change of variables and by letting $\beta=\lambda/\alpha,$
we obtain the following formula.

\begin{prop}
For $0<\alpha<1$ and $-\pi/2<\beta<\pi/2,$
$$
\int_0^{\frac\pi2-\beta}\frac{\cos^\alpha(x+\beta)}{\sin^\alpha x}dx
=\frac{\pi\sin(\frac{\pi}2-\beta)\alpha}{\sin\pi\alpha}.
$$
\end{prop}

\noindent\bigskip
{\bf Acknowledgement.}
The author is grateful to Dr.~Ikkei Hotta for his suggestions.

\def\cprime{$'$} \def\cprime{$'$} \def\cprime{$'$}
\providecommand{\bysame}{\leavevmode\hbox to3em{\hrulefill}\thinspace}
\providecommand{\MR}{\relax\ifhmode\unskip\space\fi MR }
\providecommand{\MRhref}[2]{%
  \href{http://www.ams.org/mathscinet-getitem?mr=#1}{#2}
}
\providecommand{\href}[2]{#2}

\end{document}